\documentclass{article}

\usepackage{PRIMEarxiv}
\usepackage{amsmath}
\usepackage{float}
\newtheorem{theorem}{Theorem}

\newtheorem{corollary}[theorem]{Corollary}

\usepackage{xcolor}
\usepackage[utf8]{inputenc} 
\usepackage[T1]{fontenc}    
\usepackage{hyperref}       
\usepackage{url}            
\usepackage{booktabs}       
\usepackage{amsfonts}       
\usepackage{nicefrac}       
\usepackage{microtype}      
\usepackage{lipsum}
\usepackage{fancyhdr}       
\usepackage{graphicx}       
\graphicspath{{media/}}     

\title{New Chaotic dynamics for  Yitang Zhang latest results  on Landau-Siegel zero 
\thanks{\textit{\underline{Citation}}: 
\textbf{Authors. Title. Pages.... DOI:000000/11111.}} 
}

\author{
  Zeraoulia Rafik \\
  University of batna2 Algeria \\
  Departement of mathematics \\
  Yabous, Khenchela\\
  \texttt{\{Zeraoulia .R\}r.zeraoulia@univ-batna2.dz} \\
}

\begin{document}
\maketitle

\begin{abstract}
The first part of this paper is about Consequences resulting from Yitang Zhang's latest claimed results on Landau-Siegel zero posted by some mathematicians in Mathoverflow ,For the second part we are able to derive new Chaotic dynamics for Yitang Zhang  on Landau-Siegel zero such that the behavior  of the new dynamics has been discussed ,Lyaponove Exponents has been computed and bifurcation diagram has been achieved ,The number of limit cycle and orbits are predicted.The behavior of this new dynamics roughly proves the validity of Yitang latest results .
\end{abstract}


\section{Introduction}
In mathematics, more specifically in the field of analytic number theory, a Landau–Siegel zero or simply Siegel zero (also known as exceptional zero \cite{1}, named after Edmund Landau and Carl Ludwig Siegel, is a type of potential counterexample to the generalized Riemann hypothesis, on the zeros of Dirichlet L-functions associated to quadratic number fields. Roughly speaking, these are possible zeros very near in a quantifiable sense to $s=1$ ,The way in which Siegel zeros appear in the theory of Dirichlet L-functions is as potential exceptions to the classical zero-free regions, which can only occur when the L-function is associated to a real Dirichlet character.
Recently, Yitang Zhang just gave a virtual talk about his work on Landau-Siegel zeros at Shandong University on the 5th of November's morning in China. He will also give a talk on 8th November at Peking University.He published his preprint on arXiv .
This paper \cite{2} shows that for a real primitive character $\chi$ to the modulus $D$,
$$ L(1, \chi) > c_{1}(\log D)^{-2022} $$
where $c_{1} > 0$ is an absolute, effectively computable constant.

Assuming this result is correct, what are some significant number theoretical consequences that would follow? ,Namely ,what would be the impact on PNT error estimates, arithmetic progressions, and other related problems? Does these results add somethings to the side of dynamical system and chaos theory ? .
more specifically,what is then the form of Yitang Zhang dynamics ? Is it possible to compute the entropy of this dynamics ? 

\section{First Answer by Travor Liu regarding implications on the error term of the PNT for arithmetic progressions}

It has significant implications on the error term of the PNT for arithmetic progressions.\cite{3}

\textbf{PNT and Siegel-Walfisz theorem}

Let $\psi(x;q,a)$ be the sum of $\Lambda(n)$ over $n\le x$ and $n\equiv a\pmod q$. Then the PNT states that for fixed $q$ there is

\begin{equation}\label{eq1}
\psi(x;q,a)\sim{x\over\varphi(q)}.
\end{equation}

When $q$ is not fixed, Page (1935) proved the following general result:

\begin{theorem}(page)
There exists some absolute and effective $c_0>0$ such that for all $(a,q)=1$:*

\begin{equation}\label{eq2}
\psi(x;q,a)={x\over\varphi(q)}-\color{blue}{{\chi(a)x^\beta}\over\varphi(q)\beta}+O\{xe^{-c_0\sqrt{\log x}}\}
\end{equation}
\end{theorem}
where $\chi$ denotes the exceptional character and $\beta$ denotes the Siegel zero. The blue term would be dropped if there are no exceptional characters modulo $q$.*

To unify the error terms, we require a result due to Siegel (1935):

\begin{theorem} (Siegel):
For all $\varepsilon>0$ there exists some $A_\varepsilon>0$ such that $1-\beta>A_\varepsilon q^{-\varepsilon}$.*

Plugging this result into the blue term of (\ref{eq2}), we have

\begin{equation}\label{eq3}
x^\beta\ll x e^{-A_\varepsilon q^{-\varepsilon}\log x}.
\end{equation}
If $q\le(\log x)^{2/\varepsilon}$, then the right hand side becomes $\ll xe^{-A_\varepsilon\sqrt{\log x}}$. Combining this with (\ref{eq2}) gives us the result of Walfisz (1936)
\end{theorem}

\begin{theorem} (Siegel-Walfisz):

For any $M>0$ there exists some $C_M$ such that for all $q\le(\log x)^M$ and $(a,q)=1$ there is*

\begin{equation}\label{eq4}
\psi(x;q,a)={x\over\varphi(q)}+O\{e^{-C_M\sqrt{\log x}}\},
\end{equation}
where the O-constant is absolute.
\end{theorem}
Due to the drawback in the proof of Siegel's theorem, $A_\varepsilon$ and $C_M$ are not effectively computable.

\textbf{Improvements due to Zhang}:

However, we can significantly obtain a stronger and effective improvement of Siegel-Walfisz theorem if Zhang's result is used. That is

\begin{theorem} (Zhang)
There exists $A>0$ and effective $C_1>0$ such that $L(1,\chi)>C_1(\log q)^{-A}$.
\end{theorem}
 Zhang proved this result for $A=2022$, but I choose not to plug it in for generality.
Let $\beta$ be the rightmost real zero of $L(s,\chi)$ for some real $\chi$ modulo $q$ such that $1-\beta\gg(\log q)^{-1}$. Then it follows from the mean value theorem that there exists some $1-\beta<\sigma<1$ such that $1-\beta=L(1,\chi)/L'(\sigma,\chi)$. Applying the classical bound $L'(\sigma,\chi)=O(\log^2q)$ and Zhang's result gives us the zero-free region that

$$
1-\beta>C_2(\log q)^{-A-2},
$$

where $C_2>0$ is effectively computable, which indicates that the blue term in (\ref{eq2}) is dominated by

$$
x^\beta\ll xe^{-C_2(\log x)(\log q)^{-A-2}}.
$$

If $(\log q)^{A+2}\le\sqrt{\log x}$, then the right hand side becomes $\ll xe^{-C_2\sqrt{\log x}}$, which allows Theorem 1 to be improved significantly:

\begin{theorem} 

Let $A$ be as in Theorem 4. There exists some absolute $c_0>0$ such that for all $q\le e^{(\log x)^{1/(2A+4)}}$ and $(a,q)=1$, we have*
\begin{equation}\label{eq5}
\psi(x;q,a)={x\over\varphi(q)}+O\{e^{-c_0\sqrt{\log x}}\}.
\end{equation}
\end{theorem}

\textbf{Asymptotic formulas valid for all $q\ge1$.}

Although Theorem 5 is much stronger than Theorem 3, it is difficult to compare them to Theorem 1 without the blue term, so this section is dedicated to deduce asymptotic formula valid for all $q\ge1$ and $(a,q)=1$ so that a better comparison can be made.

Since $\Lambda(n)\le\log n$, we know trivially that

\begin{equation}\label{eq6}
\psi(x;q,a)\le\sum_{\substack{n\le x\\n\equiv a(q)}}\log x\ll{x\log x\over q}.
\end{equation}

Combining this with (\ref{eq3}), we see that Theorem 3 indicates that

\begin{equation}\label{eq7}
\psi(x;q,a)={x\over\varphi(q)}+O_N\{x(\log x)^{-N}\}\quad(N>0).
\end{equation}

If the trivial upper bound is juxtaposed with (\ref{eq4}), then we see that there exists some absolute and effective $c_0>0$ such that

\begin{equation}\label{eq8}
\psi(x;q,a)={x\over\varphi(q)}+O\{xe^{-c_0(\log x)^{1/(2A+4)}}\},
\end{equation}

which has a substantially better error term than (5).

\section{Second Answer by Stopple regarding classification of discriminants of binary quadratic forms}
There will be many important consequences of Zhang's result, if correct.  One specific result is that it will reduce one of the last open problem from the era of Gauss and Euler to a finite amount of computation, namely the classification of discriminants of binary quadratic forms with one class per genus.  The congruence class of a prime number $p$ modulo $d$ determines which form of discriminant $-d<0$ represents $p$  if and only if there is one class per genus.

Such discriminants which are congruent to $0$ modulo $4$ are  Euler's *numeri idonei* or idoneal numbers.  Euler expected there would be infinitely many such discriminants.  [**Edit**: Apparently I'm mis-remembering this.  See the remarks of KConrad below.] It was Gauss who conjectured  that the only such discriminants are the 65 examples (not necessarily fundamental) known to Euler.  There are also 65 known fundamental discriminants (not necessarily even) with one class per genus.  The existence of a 66th is still an open problem.  By genus theory we know that for discriminants with one class per genus, the class group satisfies
$$
\mathcal C(-d)\cong \left(\mathbb Z/2\right)^{g-1},
$$
where $g$ is the number of prime divisors of $d$.  Obviously $d$ is bigger than the absolute value of the smallest fundamental discriminant with $g$ prime divisors,
$$
d_g\overset{\text{def.}}=3\cdot4\cdot5\cdot7\cdots p_g.
$$
From lower bounds on the size of $p_g$, the $g$-th prime, and on $\theta(x)=\sum_{p\leq x}\log(p)$, one can show that
$$
d_g>g^g.
$$
Since $
2^{g-1} \ll \sqrt{g^g}, 
$
lower bounds for the class number which we expect to be true rule out the possibility for one class per genus for large $g$.  

In 1973, Peter Weinberger showed that on GRH, no fundamental discriminant $-d<-5460$ has one class per genus, and unconditionally there is at most one more such $d$.    

In contrast, Oesterle explicitly observed that the lower bound due to Goldfeld-Gross-Zagier  is not strong enough to finish the classification of discriminants with one class per genus: 
$
\log(g^g) $ is $\ll 2^{g-1}
$.
Iwaniec and Kowalski observed that even the full strength of the Birch Swinnerton-Dyer conjecture, "the best effective lower bounds which current technology allows us to hope for"  would not suffice, as $ \log(g^g)^r$ is $\ll 2^{g-1}$  for any $r$.  In fact, the outlook is still more bleak:  Watkins observed that if the discriminant $-d$ is divisible by all the primes up to $(\log\log d)^3$ (as $d_g$ certainly is), the product over primes dividing $d$ in the Goldfeld-Gross-Zagier lower bound is so small the resulting bound is worse than the trivial bound.

If the implied constant is made explicit, Zhang's result would eliminate the possibility of one class per genus for discriminants above some bound.  For example, neglecting the constant, the lower bound is discriminants with more than 6007 prime divisors.  This works out to $d>3\cdot 10^{25734}$.

\section{Chaotic dynamics for  Yitang Zhang latest results  on Landau-Siegel zeros}
In this section we may give an explicit dynamics \cite{4} for  Yitang Zhang latest results  on Landau-Siegel zeros using both of Dirichlet definition which uses odd and even Characters $\chi$  mod $m $and Theorem 4 of Zhang also we shall discuss its behavior and we show that dynamics has a weak transition to chaos for some values and parameters of that discret new dynamics.
\begin{corollary}
For special values $s=1$, let $\chi$  be a primitive real character mod $m$, put $\xi_m=\exp \bigg(\frac{2\pi i}{m} \bigg)$, put $K=\mathbb{Q} \bigg(\sqrt{\chi(-1)m}\bigg)$, $h$ is its class number of roots of unity in it and $\epsilon$ is its fundamental unit ,Dirichelet L Function may be defined as :
\begin{equation}\label{Zhang}
L(1,\chi)=\begin{cases} 
\displaystyle\frac{2\pi h}{w\sqrt{m}} \quad  \textit{if}\ \chi(-1)=-1  \\

\displaystyle \frac{2 h \log{|\epsilon|}}{w\sqrt{m}} \quad  \textit{if}\ \chi(-1)=1 \\
\end{cases}
\end{equation}
\end{corollary}
Yitang discrete dynamics maybe derived from the definition of Dirichlet L function \cite{5} in (\ref{Zhang}) and Theorem 4 , we may start with the first  case $chi(-1)=-1$ in (\ref{Zhang}),Using both of \textbf{Theorem} 4 and (\ref{Zhang})  the 1-D discrete dynamics can be  reformulated as :

\begin{equation}\label{Dynamic1}
    x_{n+1}=\frac{\beta}{\sqrt{x_n}}+c\log(x_n)^{-\alpha},\quad \chi(-1)=-1
\end{equation}
with $x_n=m,n=0,1\cdots$ and $\beta=\frac{2\pi h}{w}$ , $c$ is the constant of Yitang defined in \textbf{Theorem} 4.

  For the second case :$\chi(-1)=1$, Yitang Dynamics can be written as :
  \begin{equation}\label{Dynamic2}
    x_{n+1}=\frac{\beta \log(|\epsilon|)}{\pi\sqrt{x_n}}+c\log(x_n)^{-\alpha},\quad \chi(-1)=1
\end{equation}
For Chaotic behavior the iteration of dynamic 1 (\ref{Dynamic1}) and dynamics 2 (\ref{Dynamic2}) 50000 times  gives the similar behavior as we will show in the following section .

\section{Analysis and discussion}

\begin{figure}[H]
    \centering
    \includegraphics{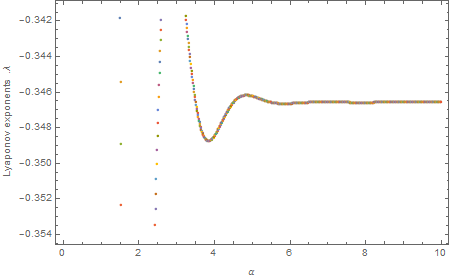}
    \caption{Lyponov exponents for the  dynamics 1  with $c\in [0.0005;0.007]$,$\alpha \in [0;10]$ }
    \label{Non chaos}
\end{figure}

For $\alpha$ in the range $[0;10]$ and  $c$ in the range $[0.0005;0.007]$ ,All Lyaponove exponents $(\lambda_1 <\lambda_2 <\lambda_3,\cdots > \lambda_{10}=-0.348 <0 $)  of the both dynamics are negative lead to the Yitang dynamics being non Chaotic ,However we noticed weak transition to Chaos in the range $\alpha \in(3;4)$ , see Figure \ref{Non chaos}

\begin{figure}[H]
    \centering
    \includegraphics{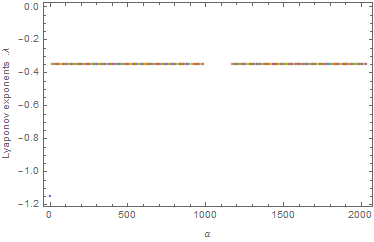}
    \caption{Lyponov exponents for the  dynamics 1  with $c\in [0.0005;0.007]$,$\alpha \in [0;2022]$ }
    \label{Non chaos2}
\end{figure}

Increasing more values of $\alpha$ up to $2022$ and $ c\in [0.0005;0.007]$, we see that the Lyponov exponents are constants approximately equals to  $-0.39$ ,see Figure \ref{Non chaos2}

\begin{figure}[H]
    \centering
    \includegraphics{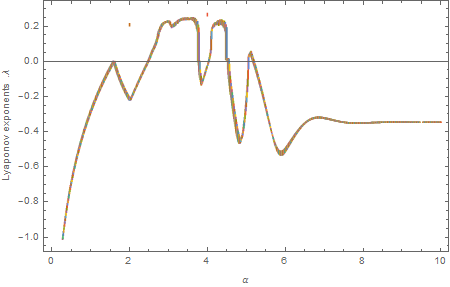}
    \caption{Lyponov exponents for the  dynamics 1  with $c\in [10.007;11.007]$,$\alpha \in (0;5)$ }
    \label{Chaos3}
\end{figure}

Transition to Chaos appear strongly  for $c\in (10.007;11.07)$ and $\alpha \in (3;5)$ such that  Lyaponove exponents values are positive $(0\leq \lambda\leq 0.21)$ ,  the max of Lyaponove exponents up to $\lambda=0.21$ see Figure \ref{Chaos3}

\begin{figure}[H]
    \centering
    \includegraphics{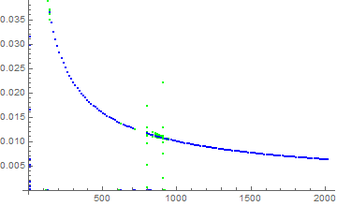}
    \caption{Bifurcation diagram for Yitang dynamics}
    \label{Chaos5}
\end{figure}

Increasing more values of $c$ one can obtain large range for transition to chaos for Yitang dynamics .The bifurcation diagram of 1.D Yitang dynamics for some values of $c$  has been plotted ,see Figure \ref{Chaos5} 

We have  discovered period doubling bifurcations by varying the parameter $\alpha$ in the chaotic range  and plotting $x_n,y_n$ versus the iterations $n$ of the map and then print below the values of the points in order to distinguish period cycles. 

\begin{figure}[H]
    \centering
    \includegraphics{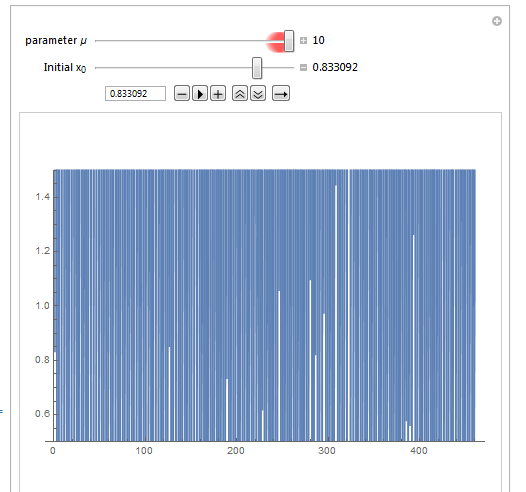}
    \caption{Yitang dynamics Fixed Points Plot versus Iterations,$c=100$,$\mu=\alpha \in (0;5)$ }
    \label{Chaos4}
\end{figure}

For one dimensional systems of the form $x_{n+1}=f(x_n)$ the period-$2$ cycles happen when the system:
$$ x_1=f(x_2) , \quad  x_2=f(x_1) $$
has a unique solution and it  seems that  $0$ is the fixed point for  Yitang dynamics (see Table1) with initial point $x_0=0.4 $. In other words, the trajetory "jumps" between $x_1$ and $x_2$ so while iterating the map they keep appearing (until a certain value of the parameter of the system). see Figure \ref{Chaos4}

\begin{table}[H]
\begin{center}
$$
\left(
\begin{array}{cc}
 \text{point} & \text{x} \\
 86 & 3.55687 \\
 87 & 0.265123 \\
 88 & 0.971062 \\
 89 & 4.05719\times 10^{76} \\
 90 & \text{2.4823161470320785$\grave{ }$*${}^{\wedge}$-39} \\
 91 & 1.00356\times 10^{19} \\
 92 & \text{1.5783353198977849$\grave{ }$*${}^{\wedge}$-10} \\
 93 & 39798.8 \\
 94 & 0.00250631 \\
 95 & 9.9874 \\
 96 & 0.158214 \\
 97 & 1.25704 \\
 98 & 1.07455\times 10^{32} \\
 99 & \text{4.823440210703131$\grave{ }$*${}^{\wedge}$-17} \\
 100 & 7.19932\times 10^7 \\
 101 & 0.0000589283 \\
\end{array}
\right)
$$
\caption{\label{demo-table} Fixed Points values of Yitang dynamics with iteration of $500$ times }
\end{center}
\end{table}

We computed the Lyapunove Exponents of Yitang Dynamics by varying alpha from 0 to 600 we got the following results in  table 2 and Plot:

Some lyapunov exponents for alpha in the range $(0,600)$ \cite{7} has been computed using mathematica as shwon in the following separated tables with plot:

\begin{table}[H]
\begin{center}
$$
\left(
\begin{array}{cc}
 \text{$\alpha$} & \text{Lyapunov exponents values} \\
 0. & -1.15226 \\
 0.01 & -1.11974 \\
 0.02 & -1.08835 \\
 0.03 & -1.058 \\
 0.04 & -1.02864 \\
 0.05 & -1.00021 \\
 0.06 & -0.972639 \\
 0.07 & -0.945889 \\
 0.08 & -0.91991 \\
 0.09 & -0.894661 \\
 0.1 & -0.870103 \\
 0.11 & -0.8462 \\
 0.12 & -0.822918 \\
 0.13 & -0.800227 \\
 0.14 & -0.778098 \\
 0.15 & -0.756505 \\
 0.16 & -0.735422 \\
 0.17 & -0.714826 \\
 0.18 & -0.694697 \\
 0.19 & -0.675013 \\
 0.2 & -0.655756 \\
 0.21 & -0.636908 \\
 0.22 & -0.618452 \\
 0.23 & -0.600373 \\
 0.24 & -0.582655 \\
 0.25 & -0.565286 \\
 0.26 & -0.54825 \\
 0.27 & -0.531537 \\
 0.28 & -0.515134 \\
 0.29 & -0.499031 \\
 0.3 & -0.483215 \\
 0.31 & -0.467679 \\
 0.32 & -0.452411 \\
 0.33 & -0.437404 \\
 0.34 & -0.422647 \\
 0.35 & -0.408134 \\
 0.36 & -0.393857 \\
 0.37 & -0.379807 \\
 0.38 & -0.365978 \\
 0.39 & -0.352363 \\
 0.4 & -0.338955 \\
 0.41 & -0.325749 \\
 0.42 & -0.312739 \\
 0.43 & -0.299918 \\
 0.44 & -0.287281 \\
 0.45 & -0.274824 \\
 0.46 & -0.262541 \\
 0.47 & -0.250427 \\
 0.48 & -0.238478 \\
 0.49 & -0.22669 \\
 0.5 & -0.215057 \\
 0.51 & -0.203577 \\
 0.52 & -0.192245 \\
 0.53 & -0.181057 \\
 0.54 & -0.17001 \\
 0.55 & -0.1591 \\
 0.56 & -0.148324 \\
 0.57 & -0.137679 \\
 0.58 & -0.127162 \\
 0.59 & -0.116769 \\
 0.6 & -0.106497 \\
 0.61 & -0.0963444 \\
 0.62 & -0.0863079 \\
 0.63 & -0.0763849 \\
 \end{array}
\right)
$$
 \caption{\label{demo-table}Table1 }
\end{center}
\end{table}
\begin{table}[H]
\begin{center}
$$
\left(
\begin{array}{cc}
 \text{$\alpha$} & \text{Lyapunov exponents values} \\
 0.64 & -0.066573 \\
 0.65 & -0.0568697 \\
 0.66 & -0.0472728 \\
 0.67 & -0.0377802 \\
 0.68 & -0.02839 \\
 0.69 & -0.0191008 \\
 0.7 & -0.00991274 \\
 0.71 & -0.000848175 \\
 0.72 & -0.016495 \\
 0.73 & -0.0342835 \\
 0.74 & -0.0518883 \\
 0.75 & -0.0693237 \\
 0.76 & -0.0866054 \\
 0.77 & -0.103749 \\
 0.78 & -0.12077 \\
 0.79 & -0.137683 \\
 0.8 & -0.154503 \\
 0.81 & -0.171244 \\
 0.82 & -0.187921 \\
 0.83 & -0.204547 \\
 0.84 & -0.221137 \\
 0.85 & -0.237705 \\
 0.86 & -0.254262 \\
 0.87 & -0.270825 \\
 0.88 & -0.287405 \\
 0.89 & -0.304016 \\
 0.9 & -0.320673 \\
 0.91 & -0.337389 \\
 0.92 & -0.354177 \\
 0.93 & -0.371052 \\
 0.94 & -0.388027 \\
 0.95 & -0.405117 \\
 0.96 & -0.422336 \\
 0.97 & -0.439701 \\
 0.98 & -0.457226 \\
 0.99 & -0.474931 \\
 \end{array}
\right)
$$
 \caption{\label{demo-table} }
\end{center}
\end{table}

\begin{table}[H]
\begin{center}
$$
\left(
\begin{array}{cc}
 \text{$\alpha$} & \text{Lyapunov exponents values} \\
 1 & -0.49284 \\
 1. & -0.49284 \\
 1.01 & -0.510944 \\
 1.02 & -0.529281 \\
 1.03 & -0.547868 \\
 1.04 & -0.566722 \\
 1.05 & -0.58586 \\
 1.06 & -0.605304 \\
 1.07 & -0.625074 \\
 1.08 & -0.645192 \\
 1.09 & -0.665681 \\
 1.1 & -0.686564 \\
 1.11 & -0.707865 \\
 1.12 & -0.729611 \\
 1.13 & -0.751828 \\
 1.14 & -0.774542 \\
 1.15 & -0.797783 \\
 1.16 & -0.821578 \\
 1.17 & -0.845957 \\
 1.18 & -0.870949 \\
 1.19 & -0.896582 \\
 1.2 & -0.922885 \\
 1.21 & -0.949883 \\
 1.22 & -0.9776 \\
 1.23 & -1.00605 \\
 1.24 & -1.03526 \\
 1.25 & -1.06522 \\
 1.26 & -1.09592 \\
 1.27 & -1.12735 \\
 1.28 & -1.15944 \\
 1.29 & -1.19213 \\
 1.3 & -1.2253 \\
 1.31 & -1.25878 \\
 1.32 & -1.29232 \\
 1.33 & -1.32562 \\
 1.34 & -1.35824 \\
 1.35 & -1.38964 \\
 1.36 & -1.41912 \\
 1.37 & -1.44583 \\
 1.38 & -1.46875 \\
 1.39 & -1.48674 \\
 1.4 & -1.49858 \\
 1.41 & -1.50307 \\
 1.42 & -1.49913 \\
 1.43 & -1.48595 \\
 1.44 & -1.46309 \\
 1.45 & -1.43055 \\
 1.46 & -1.38873 \\
 1.47 & -1.33835 \\
 1.48 & -1.28031 \\
 1.49 & -1.21558 \\
 1.5 & -1.14504 \\
 1.51 & -1.0694 \\
 1.52 & -0.989104 \\
 1.53 & -0.904293 \\
 1.54 & -0.814681 \\
 1.55 & -0.719402 \\
 1.56 & -0.616619 \\
 1.57 & -0.502533 \\
 1.58 & -0.367932 \\
 1.59 & -0.174509 \\
 \end{array}
\right)
$$
 \caption{\label{demo-table}  }
\end{center}
\end{table}
\begin{table}[H]
\begin{center}
$$
\left(
\begin{array}{cc}
 \text{$\alpha$} & \text{Lyapunov exponents values} \\
 1.6 & 0.0881717 \\
 1.61 & 0.138893 \\
 1.62 & 0.165551 \\
 1.63 & 0.202492 \\
 1.64 & 0.222027 \\
 1.65 & 0.242128 \\
 1.66 & 0.262391 \\
 1.67 & 0.262112 \\
 1.68 & 0.291839 \\
 1.69 & 0.306873 \\
 1.7 & 0.301594 \\
 1.71 & 0.332126 \\
 1.72 & 0.326654 \\
 1.73 & 0.335829 \\
 1.74 & 0.332551 \\
 1.75 & 0.309413 \\
 1.76 & 0.349316 \\
 1.77 & 0.368101 \\
 1.78 & 0.391869 \\
 1.79 & 0.3934 \\
 1.8 & 0.387423 \\
 1.81 & 0.382245 \\
 1.82 & 0.382253 \\
 1.83 & 0.38575 \\
 1.84 & 0.386156 \\
 1.85 & 0.38766 \\
 1.86 & 0.396544 \\
 1.87 & 0.405967 \\
 1.88 & 0.415302 \\
 1.89 & 0.424401 \\
 1.9 & 0.428945 \\
 1.91 & 0.431023 \\
 1.92 & 0.435248 \\
 1.93 & 0.432161 \\
 1.94 & 0.430948 \\
 1.95 & 0.428278 \\
 1.96 & 0.432447 \\
 1.97 & 0.433208 \\
 1.98 & 0.433657 \\
 1.99 & 0.437374 \\
 2 & 0.349501 \\
 2. & 0.349501 \\
 2.01 & 0.435275 \\
 2.02 & 0.435097 \\
 2.03 & 0.438616 \\
 2.04 & 0.437307 \\
 \end{array}
\right)
$$
 \caption{\label{demo-table} }
\end{center}
\end{table}

 \begin{table}[H]
\begin{center}
$$
\left(
\begin{array}{cc}
 \text{$\alpha$} & \text{Lyapunov exponents values} \\
 2.05 & 0.440017 \\
 2.06 & 0.440217 \\
 2.07 & 0.435724 \\
 2.08 & 0.441471 \\
 2.09 & 0.44623 \\
 2.1 & 0.447963 \\
 2.11 & 0.442867 \\
 2.12 & 0.44119 \\
 2.13 & 0.439389 \\
 2.14 & 0.443069 \\
 2.15 & 0.445086 \\
 2.16 & 0.442277 \\
 2.17 & 0.444541 \\
 2.18 & 0.443888 \\
 2.19 & 0.439353 \\
 2.2 & 0.435632 \\
 2.21 & 0.43046 \\
 2.22 & 0.423262 \\
 2.23 & 0.417051 \\
 2.24 & 0.41479 \\
 2.25 & 0.412113 \\
 2.26 & 0.413576 \\
 2.27 & 0.412673 \\
 2.28 & 0.413888 \\
 2.29 & 0.413859 \\
 2.3 & 0.414319 \\
 2.31 & 0.417688 \\
 2.32 & 0.42042 \\
 2.33 & 0.416738 \\
 2.34 & 0.413714 \\
 2.35 & 0.42057 \\
 2.36 & 0.41715 \\
 2.37 & 0.412618 \\
 2.38 & 0.412854 \\
 2.39 & 0.412535 \\
 2.4 & 0.414816 \\
 2.41 & 0.414834 \\
 2.42 & 0.413606 \\
 2.43 & 0.416818 \\
 2.44 & 0.41411 \\
 2.45 & 0.417853 \\
 2.46 & 0.416911 \\
 \end{array}
\right)
$$
 \caption{\label{demo-table} }
\end{center}
\end{table}
\begin{table}[H]
\begin{center}
$$
\left(
\begin{array}{cc}
 \text{$\alpha$} & \text{Lyapunov exponents values} \\
 2.47 & 0.421352 \\
 2.48 & 0.419916 \\
 2.49 & 0.424899 \\
 2.5 & 0.430376 \\
 2.51 & 0.433199 \\
 2.52 & 0.432542 \\
 2.53 & 0.428365 \\
 2.54 & 0.428214 \\
 2.55 & 0.425867 \\
 2.56 & 0.423466 \\
 2.57 & 0.421258 \\
 2.58 & 0.420319 \\
 2.59 & 0.421149 \\
 2.6 & 0.423574 \\
 2.61 & 0.421921 \\
 2.62 & 0.424841 \\
 2.63 & 0.423726 \\
 2.64 & 0.422148 \\
 2.65 & 0.420218 \\
 2.66 & 0.418454 \\
 2.67 & 0.419379 \\
 2.68 & 0.419718 \\
 2.69 & 0.41695 \\
 2.7 & 0.41443 \\
 2.71 & 0.41608 \\
 2.72 & 0.41517 \\
 2.73 & 0.410443 \\
 2.74 & 0.409386 \\
 2.75 & 0.409979 \\
 2.76 & 0.410878 \\
 2.77 & 0.407457 \\
 2.78 & 0.406276 \\
 2.79 & 0.406546 \\
 2.8 & 0.40539 \\
 2.81 & 0.413248 \\
 2.82 & 0.41515 \\
 2.83 & 0.419008 \\
 \end{array}
\right)
$$
 \caption{\label{demo-table} }
\end{center}
\end{table}

\begin{table}[H]
\begin{center}
$$
\left(
\begin{array}{cc}
 \text{$\alpha$} & \text{Lyapunov exponents values} \\
 2.84 & 0.412071 \\
 2.85 & 0.411538 \\
 2.86 & 0.399653 \\
 2.87 & 0.389743 \\
 2.88 & 0.358811 \\
 2.89 & 0.271317 \\
 2.9 & -0.0709225 \\
 2.91 & 0.204658 \\
 2.92 & 0.339491 \\
 2.93 & 0.375873 \\
 2.94 & 0.392489 \\
 2.95 & 0.405362 \\
 2.96 & 0.415369 \\
 2.97 & 0.420442 \\
 2.98 & 0.423161 \\
 2.99 & 0.427594 \\
 3 & 0.4321 \\
 3. & 0.4321 \\
 3.01 & 0.435507 \\
 3.02 & 0.441413 \\
 3.03 & 0.444938 \\
 3.04 & 0.445524 \\
 3.05 & 0.448535 \\
 3.06 & 0.45017 \\
 3.07 & 0.452593 \\
 3.08 & 0.456965 \\
 3.09 & 0.454544 \\
 3.1 & 0.457091 \\
 3.11 & 0.453177 \\
 3.12 & 0.456127 \\
 3.13 & 0.45441 \\
 3.14 & 0.45515 \\
 3.15 & 0.456125 \\
 3.16 & 0.452111 \\
 3.17 & 0.4553 \\
 3.18 & 0.458551 \\
 3.19 & 0.457554 \\
 3.2 & 0.455247 \\
 3.21 & 0.45498 \\
 3.22 & 0.454235 \\
 3.23 & 0.455603 \\
 3.24 & 0.462065 \\
 3.25 & 0.469304 \\
 3.26 & 0.465601 \\
 3.27 & 0.466889 \\
 3.28 & 0.465339 \\
 3.29 & 0.465456 \\
 3.3 & 0.458958 \\
 3.31 & 0.459289 \\
 3.32 & 0.461778 \\
 3.33 & 0.463729 \\
 3.34 & 0.460103 \\
 3.35 & 0.4594 \\
 3.36 & 0.458898 \\
 3.37 & 0.462396 \\
 3.38 & 0.465787 \\
 3.39 & 0.459454 \\
 3.4 & 0.461186 \\
 3.41 & 0.461277 \\
 3.42 & 0.462188 \\
 3.43 & 0.458111 \\
 3.44 & 0.456439 \\
 3.45 & 0.426665 \\
 3.46 & 0.451274 \\
  \end{array}
\right)
$$
 \caption{\label{demo-table}  }
\end{center}
\end{table}

\begin{table}[H]
\begin{center}
$$
\left(
\begin{array}{cc}
 \text{$\alpha$} & \text{Lyapunov exponents values} \\
 3.47 & 0.464362 \\
 3.48 & 0.466476 \\
 3.49 & 0.469612 \\
 3.5 & 0.472629 \\
 3.51 & 0.479088 \\
 3.52 & 0.479697 \\
 3.53 & 0.485185 \\
 3.54 & 0.492687 \\
 3.55 & 0.500078 \\
 3.56 & 0.502497 \\
 3.57 & 0.502803 \\
 3.58 & 0.505774 \\
 3.59 & 0.508774 \\
 3.6 & 0.508004 \\
 3.61 & 0.515266 \\
 3.62 & 0.517609 \\
 3.63 & 0.517886 \\
 3.64 & 0.523134 \\
 3.65 & 0.527244 \\
 3.66 & 0.533924 \\
 3.67 & 0.539092 \\
 3.68 & 0.54134 \\
 3.69 & 0.538396 \\
 3.7 & 0.544749 \\
 3.71 & 0.547773 \\
 3.72 & 0.549363 \\
 3.73 & 0.551562 \\
 3.74 & 0.553626 \\
 3.75 & 0.556052 \\
 3.76 & 0.557295 \\
 3.77 & 0.56001 \\
 3.78 & 0.567793 \\
 3.79 & 0.568275 \\
 3.8 & 0.569386 \\
 3.81 & 0.571111 \\
 3.82 & 0.572617 \\
 3.83 & 0.575607 \\
 3.84 & 0.578082 \\
 3.85 & 0.58003 \\
 3.86 & 0.577699 \\
 3.87 & 0.578298 \\
 3.88 & 0.579757 \\
 3.89 & 0.579385 \\
 3.9 & 0.580589 \\
 3.91 & 0.580969 \\
 3.92 & 0.582456 \\
 3.93 & 0.576281 \\
 3.94 & 0.576627 \\
  \end{array}
\right)
$$
 \caption{\label{demo-table}}
\end{center}
\end{table}
\begin{table}[H]
\begin{center}
$$
\left(
\begin{array}{cc}
 \text{$\alpha$} & \text{Lyapunov exponents values} \\
 3.95 & 0.573432 \\
 3.96 & 0.574313 \\
 3.97 & 0.57521 \\
 3.98 & 0.57762 \\
 3.99 & 0.576215 \\
 4 & 0.507632 \\
 4. & 0.507632 \\
 4.01 & 0.576928 \\
 4.02 & 0.579918 \\
 4.03 & 0.583196 \\
 4.04 & 0.580391 \\
 4.05 & 0.580994 \\
 4.06 & 0.582307 \\
 4.07 & 0.582366 \\
 4.08 & 0.583606 \\
 4.09 & 0.585411 \\
 4.1 & 0.584807 \\
 4.11 & 0.584631 \\
 4.12 & 0.58586 \\
 4.13 & 0.586655 \\
 4.14 & 0.590236 \\
 4.15 & 0.591261 \\
 4.16 & 0.591472 \\
 4.17 & 0.591893 \\
 4.18 & 0.589417 \\
 \end{array}
\right)
$$
 \caption{\label{demo-table} }
\end{center}
\end{table}

\begin{table}[H]
\begin{center}
$$
\left(
\begin{array}{cc}
 \text{$\alpha$} & \text{Lyapunov exponents values} \\
 4.19 & 0.58858 \\
 4.2 & 0.587385 \\
 4.21 & 0.586986 \\
 4.22 & 0.586748 \\
 4.23 & 0.588275 \\
 4.24 & 0.584422 \\
 4.25 & 0.585339 \\
 4.26 & 0.583006 \\
 4.27 & 0.582097 \\
 4.28 & 0.581886 \\
 4.29 & 0.582088 \\
 4.3 & 0.581241 \\
 4.31 & 0.581622 \\
 4.32 & 0.581417 \\
 4.33 & 0.580996 \\
 4.34 & 0.579168 \\
 4.35 & 0.581486 \\
 4.36 & 0.580001 \\
 4.37 & 0.575849 \\
 4.38 & 0.569656 \\
 4.39 & 0.569377 \\
 4.4 & 0.569891 \\
 4.41 & 0.567187 \\
 4.42 & 0.570189 \\
 4.43 & 0.570313 \\
 4.44 & 0.571483 \\
 4.45 & 0.571226 \\
 4.46 & 0.574592 \\
 4.47 & 0.573314 \\
 4.48 & 0.571783 \\
 4.49 & 0.571822 \\
 4.5 & 0.569118 \\
 4.51 & 0.568564 \\
 4.52 & 0.567916 \\
 4.53 & 0.568978 \\
 4.54 & 0.570392 \\
 4.55 & 0.572938 \\
 4.56 & 0.571468 \\
 4.57 & 0.572372 \\
 4.58 & 0.571267 \\
 4.59 & 0.572128 \\
 4.6 & 0.569822 \\
 4.61 & 0.571771 \\
 4.62 & 0.574242 \\
 4.63 & 0.576563 \\
 4.64 & 0.578604 \\
 4.65 & 0.582042 \\
 4.66 & 0.581638 \\
 4.67 & 0.58133 \\
 4.68 & 0.583655 \\
 4.69 & 0.584481 \\
 4.7 & 0.586956 \\
 4.71 & 0.58891 \\
 4.72 & 0.58939 \\
 4.73 & 0.589091 \\
 4.74 & 0.590923 \\
 4.75 & 0.591047 \\
 4.76 & 0.588998 \\
 4.77 & 0.590917 \\
 4.78 & 0.588484 \\
 4.79 & 0.589111 \\
  \end{array}
\right)
$$
 \caption{\label{demo-table}  }
\end{center}
\end{table}

 \begin{table}[H]
\begin{center}
$$
\left(
\begin{array}{cc}
 \text{$\alpha$} & \text{Lyapunov exponents values} \\
 4.8 & 0.590603 \\
 4.81 & 0.590248 \\
 4.82 & 0.587849 \\
 4.83 & 0.586583 \\
 4.84 & 0.584447 \\
 4.85 & 0.579494 \\
 4.86 & 0.579748 \\
 4.87 & 0.57695 \\
 4.88 & 0.5767 \\
 4.89 & 0.571527 \\
 4.9 & 0.575459 \\
 4.91 & 0.569518 \\
 4.92 & 0.57061 \\
 4.93 & 0.569778 \\
 4.94 & 0.570774 \\
 4.95 & 0.571756 \\
 4.96 & 0.572801 \\
 4.97 & 0.570044 \\
 4.98 & 0.568464 \\
 4.99 & 0.567003 \\
 5 & 0.564626 \\
 6 & 0.510146 \\
 7 & 0.556004 \\
 8 & -0.0719482 \\
 9 & 0.534278 \\
 10 & 0.522649 \\
 11 & 0.517388 \\
 12 & 0.472097 \\
 13 & 0.492431 \\
 14 & 0.446779 \\
 15 & 0.460998 \\
 16 & 0.459563 \\
 17 & 0.458949 \\
 18 & 0.447865 \\
 19 & 0.405761 \\
 20 & -0.0390125 \\
 21 & 0.388358 \\
 22 & 0.365367 \\
 23 & 0.379947 \\
 24 & 0.337853 \\
 25 & 0.385089 \\
 26 & 0.377315 \\
 27 & 0.393322 \\
 28 & 0.384273 \\
 29 & 0.404729 \\
 30 & 0.394806 \\
 31 & 0.392084 \\
 32 & 0.35981 \\
 33 & 0.361309 \\
 34 & 0.335672 \\
 35 & 0.325758 \\
 36 & 0.332982 \\
 37 & 0.326919 \\
 38 & 0.334096 \\
 39 & 0.325021 \\
 40 & 0.31735 \\
 41 & 0.314913 \\
 42 & 0.315999 \\
 43 & 0.304623 \\
 44 & 0.297582 \\
 45 & 0.310218 \\
  \end{array}
\right)
$$
 \caption{\label{demo-table} }
\end{center}
\end{table}

\begin{table}[H]
\begin{center}
$$
\left(
\begin{array}{cc}
 \text{$\alpha$} & \text{Lyapunov exponents values} \\
 46 & 0.31658 \\
 47 & 0.312856 \\
 48 & 0.304048 \\
 49 & 0.319614 \\
 50 & 0.318441 \\
 51 & 0.326547 \\
 52 & 0.326991 \\
 53 & 0.333725 \\
 54 & 0.336305 \\
 55 & 0.325113 \\
 56 & 0.319471 \\
 57 & 0.284135 \\
 58 & -0.0551551 \\
 59 & 0.250207 \\
 60 & 0.052976 \\
 61 & 0.268896 \\
 62 & 0.243243 \\
 63 & 0.277856 \\
 64 & 0.276556 \\
 65 & 0.272771 \\
 66 & 0.261744 \\
 67 & 0.268161 \\
 68 & 0.268744 \\
 69 & 0.261705 \\
 70 & 0.250056 \\
 71 & 0.261992 \\
 72 & 0.262499 \\
 73 & 0.250503 \\
 74 & 0.252945 \\
 75 & 0.249279 \\
 76 & 0.248871 \\
 77 & 0.246984 \\
 78 & 0.253233 \\
 79 & 0.243082 \\
 80 & 0.222819 \\
 81 & 0.242301 \\
 82 & 0.235616 \\
 83 & 0.244237 \\
 84 & 0.261487 \\
 85 & 0.25179 \\
 86 & 0.25768 \\
 87 & 0.253262 \\
 88 & 0.258194 \\
 89 & 0.253025 \\
 90 & 0.244901 \\
 91 & 0.252182 \\
 92 & 0.261088 \\
 93 & 0.26153 \\
 94 & 0.283992 \\
 95 & 0.268048 \\
 96 & 0.278632 \\
 97 & 0.269725 \\
 98 & 0.265437 \\
 99 & 0.278103 \\
 \end{array}
\right)
$$
 \caption{\label{demo-table}}
\end{center}
\end{table}
  
 \begin{table}[H]
\begin{center}
$$
\left(
\begin{array}{cc}
 \text{$\alpha$} & \text{Lyapunov exponents values} \\
 100 & 0.294507 \\
 101 & 0.280263 \\
 102 & 0.338833 \\
 103 & 0.282768 \\
 104 & 0.351269 \\
 105 & 0.271341 \\
 106 & 0.28702 \\
 107 & 0.258217 \\
 108 & 0.298877 \\
 109 & 0.237881 \\
 110 & 0.267863 \\
 111 & 0.243872 \\
 112 & 0.312845 \\
 113 & 0.243733 \\
 114 & 0.341542 \\
 115 & 0.248787 \\
 116 & 0.317934 \\
 117 & 0.247616 \\
 118 & 0.361362 \\
 119 & 0.245811 \\
 120 & 0.270793 \\
 121 & 0.238476 \\
 122 & 0.264539 \\
 123 & 0.240243 \\
 124 & 0.241082 \\
 125 & 0.237228 \\
 126 & 0.306439 \\
 127 & 0.240514 \\
 128 & 0.253423 \\
 129 & 0.233118 \\
 130 & 0.363916 \\
 131 & 0.235455 \\
 132 & 0.236588 \\
 133 & 0.232666 \\
 134 & 0.239603 \\
 135 & 0.228082 \\
 136 & -0.000682791 \\
 137 & 0.221024 \\
 138 & 0.248629 \\
 139 & 0.233006 \\
 140 & 0.249864 \\
 141 & 0.218172 \\
 142 & 0.246957 \\
 143 & 0.222976 \\
 144 & 0.330647 \\
 145 & 0.222045 \\
 146 & 0.227534 \\
 147 & 0.209626 \\
 148 & 0.224004 \\
 149 & 0.209111 \\
 150 & 0.238807 \\
 151 & 0.207961 \\
 152 & 0.23251 \\
 153 & 0.204123 \\
 154 & 0.258405 \\
 155 & 0.203154 \\
 156 & 0.261367 \\
 \end{array}
\right)
$$
 \caption{\label{demo-table}  }
\end{center}
\end{table}
\begin{table}[H]
\begin{center}
$$
\left(
\begin{array}{cc}
 \text{$\alpha$} & \text{Lyapunov exponents values} \\
 157 & 0.200024 \\
 158 & 0.276927 \\
 159 & 0.212927 \\
 160 & 0.227649 \\
 161 & 0.205282 \\
 162 & 0.244713 \\
 163 & 0.215555 \\
 164 & 0.224634 \\
 165 & 0.210837 \\
 166 & 0.218134 \\
 167 & 0.213816 \\
 168 & 0.25477 \\
 169 & 0.209253 \\
 170 & 0.338271 \\
 171 & 0.219416 \\
 172 & 0.318408 \\
 173 & 0.209499 \\
 174 & 0.274409 \\
 175 & 0.207832 \\
 176 & 0.278909 \\
 177 & 0.219375 \\
 178 & 0.219368 \\
 179 & 0.21559 \\
 180 & 0.266607 \\
 181 & 0.217853 \\
 182 & 0.228786 \\
 183 & 0.225598 \\
 184 & 0.23511 \\
 185 & 0.231376 \\
 186 & 0.256866 \\
 187 & 0.245295 \\
 188 & 0.230966 \\
 189 & 0.244057 \\
 \end{array}
\right)
$$
 \caption{\label{demo-table} }
\end{center}
\end{table}

\begin{table}[H]
\begin{center}
$$
\left(
\begin{array}{cc}
 \text{$\alpha$} & \text{Lyapunov exponents values} \\
 190 & 0.267777 \\
 191 & 0.250117 \\
 192 & 0.302975 \\
 193 & 0.255108 \\
 194 & 0.285148 \\
 195 & 0.257763 \\
 196 & 0.28807 \\
 197 & 0.243462 \\
 198 & 0.318595 \\
 199 & 0.203865 \\
 200 & -0.000811814 \\
 201 & 0.0151214 \\
 202 & -0.000277031 \\
 203 & 0.200983 \\
 204 & 0.0421974 \\
 205 & 0.246805 \\
 206 & 0.238059 \\
 207 & 0.22248 \\
 208 & 0.259055 \\
 209 & 0.240073 \\
 210 & 0.244049 \\
 211 & 0.235992 \\
 212 & 0.214536 \\
 213 & 0.23179 \\
 214 & 0.227318 \\
 215 & 0.224653 \\
 216 & 0.265822 \\
 217 & 0.213538 \\
 218 & 0.222347 \\
 219 & 0.238976 \\
 220 & 0.261832 \\
 221 & 0.234598 \\
 222 & 0.266934 \\
 223 & 0.229075 \\
 224 & 0.204892 \\
 225 & 0.243175 \\
 226 & 0.256176 \\
 227 & 0.216407 \\
 228 & 0.311056 \\
 229 & 0.212253 \\
 230 & 0.245005 \\
 231 & 0.212124 \\
 232 & 0.240131 \\
 233 & 0.222259 \\
 234 & 0.254617 \\
 235 & 0.199404 \\
 236 & 0.212442 \\
 237 & 0.220351 \\
 238 & 0.260493 \\
 239 & 0.221093 \\
 240 & 0.219656 \\
 241 & 0.205173 \\
 \end{array}
\right)
$$
 \caption{\label{demo-table}}
\end{center}
\end{table}
 
 \begin{table}[H]
\begin{center}
$$
\left(
\begin{array}{cc}
 \text{$\alpha$} & \text{Lyapunov exponents values} \\
 242 & 0.2216 \\
 243 & 0.192342 \\
 244 & 0.250546 \\
 245 & 0.204687 \\
 246 & 0.262024 \\
 247 & 0.191431 \\
 248 & 0.268765 \\
 249 & 0.204364 \\
 250 & 0.239521 \\
 251 & 0.198227 \\
 252 & 0.24595 \\
 253 & 0.207629 \\
 254 & 0.213018 \\
 255 & 0.203958 \\
 256 & 0.230494 \\
 257 & 0.202481 \\
 258 & 0.235047 \\
 259 & 0.196655 \\
 260 & 0.206387 \\
 261 & 0.236002 \\
 262 & 0.26797 \\
 263 & 0.213241 \\
 264 & 0.235129 \\
 265 & 0.226122 \\
 266 & 0.220672 \\
 267 & 0.182704 \\
 268 & 0.21489 \\
 269 & 0.229958 \\
 270 & 0.222978 \\
 271 & 0.198288 \\
 272 & 0.233587 \\
 273 & 0.223525 \\
 274 & 0.198087 \\
 275 & 0.203532 \\
 276 & 0.253099 \\
 277 & 0.186744 \\
 278 & 0.217099 \\
 279 & 0.248554 \\
 280 & 0.201609 \\
 281 & 0.213417 \\
 282 & 0.215235 \\
 283 & 0.199976 \\
 284 & 0.286658 \\
 285 & 0.216406 \\
 286 & 0.231279 \\
 287 & 0.181706 \\
 288 & 0.219251 \\
 289 & 0.201476 \\
 290 & 0.292975 \\
 291 & 0.203687 \\
 292 & 0.197124 \\
 293 & 0.209663 \\
 294 & 0.194105 \\
  \end{array}
\right)
$$
 \caption{\label{demo-table} }
\end{center}
\end{table}
\begin{table}[H]
\begin{center}
$$
\left(
\begin{array}{cc}
 \text{$\alpha$} & \text{Lyapunov exponents values} \\
 295 & 0.174838 \\
 296 & 0.218176 \\
 297 & 0.20286 \\
 298 & 0.229976 \\
 299 & 0.200354 \\
 300 & 0.232252 \\
 301 & 0.193725 \\
 302 & 0.221419 \\
 303 & 0.188222 \\
 304 & 0.204547 \\
 305 & 0.200112 \\
 306 & 0.218766 \\
 307 & 0.237877 \\
 308 & 0.220502 \\
 309 & 0.182346 \\
 310 & 0.210023 \\
 311 & 0.239338 \\
 312 & 0.188054 \\
 313 & 0.221709 \\
 314 & 0.209987 \\
 315 & 0.177275 \\
 316 & 0.255036 \\
 317 & 0.224455 \\
 318 & 0.238405 \\
 319 & 0.185702 \\
 320 & 0.231805 \\
 321 & 0.211784 \\
 322 & 0.210903 \\
 323 & 0.248652 \\
 324 & 0.274243 \\
 325 & 0.238526 \\
 326 & 0.227284 \\
 327 & 0.211009 \\
 328 & 0.222291 \\
 329 & 0.211565 \\
 330 & 0.230572 \\
 \end{array}
\right)
$$
 \caption{\label{demo-table}}
\end{center}
\end{table}

\begin{table}[H]
\begin{center}
$$
\left(
\begin{array}{cc}
 \text{$\alpha$} & \text{Lyapunov exponents values} \\
 331 & 0.20134 \\
 332 & 0.219904 \\
 333 & 0.213297 \\
 334 & 0.219925 \\
 335 & 0.209362 \\
 336 & 0.266719 \\
 337 & 0.215589 \\
 338 & 0.249748 \\
 339 & 0.261866 \\
 340 & 0.236032 \\
 341 & 0.247475 \\
 342 & 0.215065 \\
 343 & 0.207613 \\
 344 & 0.229582 \\
 345 & 0.21395 \\
 346 & 0.227397 \\
 347 & 0.227974 \\
 348 & 0.260846 \\
 349 & 0.194388 \\
 350 & 0.183481 \\
 351 & 0.201502 \\
 352 & 0.211913 \\
 353 & 0.229533 \\
 354 & 0.228183 \\
 355 & 0.190445 \\
 356 & 0.238333 \\
 357 & 0.219304 \\
 358 & 0.251103 \\
 359 & 0.241116 \\
 360 & 0.219417 \\
 361 & 0.202869 \\
 362 & 0.224552 \\
  \end{array}
\right)
$$
 \caption{\label{demo-table} }
\end{center}
\end{table}
 
 \begin{table}[H]
\begin{center}
$$
\left(
\begin{array}{cc}
 \text{$\alpha$} & \text{Lyapunov exponents values} \\
 363 & 0.230416 \\
 364 & 0.240793 \\
 365 & 0.223057 \\
 366 & 0.260618 \\
 367 & 0.226144 \\
 368 & 0.291342 \\
 369 & 0.230966 \\
 370 & 0.267265 \\
 371 & 0.290253 \\
 372 & 0.249541 \\
 373 & 0.240764 \\
 374 & 0.272625 \\
 375 & 0.219788 \\
 376 & 0.196459 \\
 377 & 0.286189 \\
 378 & 0.335753 \\
 379 & 0.304186 \\
 380 & 0.304906 \\
 381 & 0.300879 \\
 382 & 0.328899 \\
 383 & 0.313488 \\
 384 & 0.275138 \\
 385 & 0.320551 \\
 386 & 0.356407 \\
 387 & 0.298877 \\
 388 & 0.112348 \\
 389 & 0.289324 \\
 390 & 0.0848054 \\
 391 & 0.302225 \\
 392 & 0.207087 \\
 393 & 0.38845 \\
 394 & 3.251 \\
 395 & 0.416102 \\
 396 & 0.353294 \\
 397 & 0.302942 \\
 398 & 0.315576 \\
 399 & 0.335486 \\
 400 & 0.331842 \\
 401 & 0.286329 \\
 402 & 0.285197 \\
 403 & 0.34339 \\
 404 & 0.295304 \\
 405 & 0.299275 \\
 406 & 0.768472 \\
 407 & 0.369837 \\
 \end{array}
\right)
$$
 \caption{\label{demo-table}  }
\end{center}
\end{table}

\begin{table}[H]
\begin{center}
$$
\left(
\begin{array}{cc}
 \text{$\alpha$} & \text{Lyapunov exponents values} \\
 408 & 0.349883 \\
 409 & 0.31186 \\
 410 & 0.777299 \\
 411 & 0.312361 \\
 412 & 0.331389 \\
 413 & 0.284097 \\
 414 & 0.28867 \\
 415 & 0.310068 \\
 416 & 0.304586 \\
 417 & 0.271222 \\
 418 & 0.221264 \\
 419 & 0.324382 \\
 420 & 0.283068 \\
 421 & 0.253938 \\
 422 & 0.316774 \\
 423 & 0.294786 \\
 424 & 0.296038 \\
 425 & 0.262009 \\
 426 & 0.871476 \\
 427 & 0.333812 \\
 428 & 0.221338 \\
 429 & 0.316276 \\
 430 & 0.299743 \\
 431 & 0.34594 \\
 432 & 0.279604 \\
 433 & 0.271286 \\
 434 & 0.270754 \\
 435 & 0.258876 \\
 436 & 0.285789 \\
 437 & 0.243558 \\
 438 & 0.280651 \\
 439 & 0.273258 \\
 440 & 0.365818 \\
 441 & 0.262362 \\
 442 & 1.08994 \\
 443 & 0.267152 \\
 444 & 0.265681 \\
 445 & 0.272178 \\
 \end{array}
\right)
$$
 \caption{\label{demo-table} }
\end{center}
\end{table}

 \begin{table}[H]
\begin{center}
$$
\left(
\begin{array}{cc}
 \text{$\alpha$} & \text{Lyapunov exponents values} \\
 446 & 0.147717 \\
 447 & 0.266351 \\
 448 & 0.310172 \\
 449 & 0.290583 \\
 450 & 0.183513 \\
 451 & 0.275272 \\
 452 & 0.325405 \\
 453 & 0.276304 \\
 454 & 0.185019 \\
 455 & 0.361066 \\
 456 & 0.34653 \\
 457 & 0.283958 \\
 458 & 0.275742 \\
 459 & 0.323037 \\
 460 & 0.326334 \\
 461 & 0.289564 \\
 462 & 0.343046 \\
 463 & 0.322932 \\
 464 & 0.59162 \\
 465 & 0.26335 \\
 466 & 0.293619 \\
 467 & 0.267087 \\
 468 & 0.31462 \\
 469 & 0.351324 \\
 470 & 0.411217 \\
 471 & 0.298389 \\
 472 & 0.721084 \\
 473 & 0.280495 \\
 474 & 0.308878 \\
 475 & 0.297671 \\
 476 & 0.272683 \\
 477 & 0.295723 \\
 478 & 0.293 \\
 479 & 0.334672 \\
 480 & 0.397376 \\
 \end{array}
\right)
$$
 \caption{\label{demo-table} }
\end{center}
\end{table}

 \begin{table}[H]
\begin{center}
$$
\left(
\begin{array}{cc}
 \text{$\alpha$} & \text{Lyapunov exponents values} \\
 481 & 0.317711 \\
 482 & 0.238976 \\
 483 & 0.308637 \\
 484 & 0.296584 \\
 485 & 0.263766 \\
 486 & 0.248706 \\
 487 & 0.254955 \\
 488 & 0.250237 \\
 489 & 0.348082 \\
 490 & 0.167979 \\
 491 & 0.304352 \\
 492 & 0.448908 \\
 493 & 0.291155 \\
 494 & 0.313148 \\
 495 & 0.265916 \\
 496 & 0.19097 \\
 497 & 0.287631 \\
 498 & 0.369616 \\
 499 & 0.233645 \\
 500 & 0.345756 \\
 501 & 0.369985 \\
 502 & 0.30685 \\
 503 & 0.255885 \\
 504 & 0.569668 \\
 505 & 0.289266 \\
 506 & 0.262803 \\
 507 & 0.31031 \\
 508 & 0.338266 \\
 509 & 0.351184 \\
 510 & 0.390376 \\
 511 & 0.287916 \\
 512 & 0.27619 \\
 513 & 0.273717 \\
 514 & 0.551982 \\
 515 & 0.306474 \\
 516 & 0.316648 \\
 517 & 0.275735 \\
 518 & 0.377191 \\
 519 & 0.301323 \\
 520 & 0.295325 \\
 521 & 0.329195 \\
 \end{array}
\right)
$$
 \caption{\label{demo-table}  }
\end{center}
\end{table}
 
 \begin{table}[H]
\begin{center}
$$
\left(
\begin{array}{cc}
 \text{$\alpha$} & \text{Lyapunov exponents values} \\
 522 & 0.301517 \\
 523 & 0.355049 \\
 524 & 0.45739 \\
 525 & 0.267099 \\
 526 & 0.336464 \\
 527 & 0.232015 \\
 528 & 0.302149 \\
 529 & 0.261932 \\
 530 & 0.293415 \\
 531 & 0.365807 \\
 532 & 0.153578 \\
 533 & 0.316334 \\
 534 & 0.340326 \\
 535 & 0.274738 \\
 536 & 0.368267 \\
 537 & 0.320209 \\
 538 & 0.60179 \\
 539 & 0.264732 \\
 540 & 0.282688 \\
 541 & 0.220575 \\
 542 & 0.403453 \\
 543 & 0.304669 \\
 544 & 0.326846 \\
 545 & 0.305038 \\
 546 & 0.390608 \\
 547 & 0.313079 \\
 548 & 0.270487 \\
 549 & 0.349348 \\
 550 & 0.278621 \\
 551 & 0.29754 \\
 552 & 0.349668 \\
 553 & 0.292467 \\
 554 & 0.258608 \\
 555 & 0.305177 \\
 556 & 0.153286 \\
 557 & 0.154885 \\
 558 & 0.50595 \\
 559 & 0.315583 \\
 560 & 0.146772 \\
 561 & 0.299072 \\
  \end{array}
\right)
$$
 \caption{\label{demo-table}  }
\end{center}
\end{table}

\begin{table}[H]
\begin{center}
$$
\left(
\begin{array}{cc}
 \text{$\alpha$} & \text{Lyapunov exponents values} \\
 562 & 0.433062 \\
 563 & 0.253575 \\
 564 & 0.197195 \\
 565 & 0.295641 \\
 566 & 0.348164 \\
 567 & 0.307269 \\
 568 & 0.309574 \\
 569 & 0.364853 \\
 570 & 0.256416 \\
 571 & 0.314278 \\
 572 & 0.138215 \\
 573 & 0.33564 \\
 574 & 0.301685 \\
 575 & 0.291345 \\
 576 & 0.26482 \\
 577 & 0.380644 \\
 578 & 0.152456 \\
 579 & 0.317665 \\
 580 & 0.397127 \\
 581 & 0.327612 \\
 582 & 0.244467 \\
 583 & 0.305 \\
 584 & 0.247417 \\
 585 & 0.301816 \\
 586 & 0.227344 \\
 587 & 0.296792 \\
 588 & 0.721956 \\
 589 & 0.337911 \\
 590 & 0.257894 \\
 591 & 0.28919 \\
 592 & 0.298331 \\
 593 & 0.268303 \\
 594 & 0.337063 \\
 595 & 0.319528 \\
 596 & 0.43679 \\
 597 & 0.306326 \\
 598 & 0.477708 \\
 599 & 0.28166 \\
 600 & 0.348162 \\
\end{array}
\right)
$$
 \caption{\label{demo-table} }
\end{center}
\end{table}

The above listed lyapunov exponents are interpreted in the following plot (see Figure \ref{Non chaos.y}:
\begin{figure}[H]
    \centering
    \includegraphics{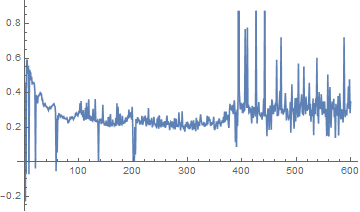}
    \caption{Lyponov exponents for the  dynamics 1 ,$\alpha \in [0;600]$ }
    \label{Non chaos.y}
\end{figure}

Yitang in \cite{6}  established a relationship between the existence of the Landau-Siegel zero of $L(s,\chi)$ and the distribution of zeros of the Dirichlet $L$-function $L(s,\psi)$, with $\psi$ belonging to a set $\Psi$ of primitive characters, in a region $\Omega$. He showed  that if the Landau-Siegel zero exists (equivalently, $L(1,\chi)$ is small),Now assume that our dynamics is reformulated correctly ,our computations showed that dynamics is chaotic  ,now we may consider the analytical solutions of Discrete dynamic 1 (Yitang Dynamics) are exactly  the zeros of $L(1,\chi)$ , The above numerical evidence showed that the Yitang discrete dynamics has a fixed point $\sigma$ close to $0$ but not exactly $0$ for $L(1,\chi)$ ($\alpha$ large ), where the system is non chaotic for $\alpha$ large enough , Assume that  $ 1-\beta<\sigma<1$ ,follows from the mean value theorem (since the dynamics is non linear) ,$1-\beta=L(1,\chi)/L'(\sigma,\chi)$. Applying the classical bound $L'(\sigma,\chi)=O(\log^2q)$ and Zhang's result (Theorem 4) gives us the zero-free region that

$$
1-\beta>C_2(\log q)^{-A-2},
$$
 and mayeb this will influence the zeros of other $L$-functions, forcing them to lie on the critical line (which is expected) and to be very regularly spaced because for large $\alpha$  the dynamics is non chaotic lead to the zeros of $L(s,\chi)$ uniforme distributed .

\section{Conclusion}
We have discussed  a new discrete dynamical system that affirms the results obtained by Yitang  such that the computation of Lyapunove exponents of the  new dynamics showed us that no transition to chaos for large values of $  \alpha$ which means no chaos for small $L(1,\chi)$ and there is transition to chaos for small  $\alpha$ in small range  .The behavior of this new dynamics roughly proves the validity of Yitang latest results on  the Landau-Siegel zero .

\section*{Acknowledgments}
This research paper is a voluntary work for researchers who do not use the Math- Overflow website, as it summarized some of the results and showed the significant positive progress in the case of the validity of the results obtained by Yitang Zhang about the Landau zero, while this paper highlighted a new discrete system whose chaotic behavior roughly proves of Yitang latest results

\bibliographystyle{unsrt}

\end{document}